\documentclass[11pt]{article}

\usepackage{amscd,amsmath, amssymb, fancyhdr, graphicx, url}

\numberwithin{equation}{section}

\newcommand{\version}{version 1.0,\ \ Nov. 7, 2015}

\renewcommand{\leftmark}{{\scriptsize Hyperbolic geometry of the ample cone}}

\makeatletter
\def\x@arrow{\DOTSB\Relbar}
\def\xlongrightarrowfill@{\arrowfill@\relbar\relbar\longrightarrow}
\newcommand{\xlongrightarrow}[2][]{%
        \ext@arrow 0099\xlongrightarrowfill@{#1}{#2}}
\makeatother

\def\eqref#1{(\ref{#1})}

\newcommand{\Z}{{\Bbb Z}}
\def\C{{\Bbb C}}
\def\P{{\Bbb P}}

\newcommand{\R}{{\Bbb R}}
\newcommand{\Q}{{\Bbb Q}}
\renewcommand{\H}{{\Bbb H}}
\newcommand{\6}{\partial}
\def\1{\sqrt{-1}\:}

\newcommand{\cntrct}                
{\hspace{2pt}\raisebox{1pt}{\text{$\lrcorner$}}\hspace{2pt}}

\renewcommand{\phi}{\varphi}
\renewcommand{\epsilon}{\varepsilon}
\renewcommand{\geq}{\geqslant}

\newcommand{\Teich}{\operatorname{\sf Teich}}

\newcommand{\Comp}{\operatorname{\sf Comp}}

\newcommand{\Kah}{\operatorname{Kah}}
\newcommand{\Amp}{\operatorname{Amp}}

\newcommand{\Pos}{\operatorname{Pos}}

\newcommand{\Aut}{\operatorname{Aut}}

\newcommand{\Diff}{\operatorname{\sf Diff}}

\newcommand{\St}{\operatorname{St}}

\newcommand{\Hdg}{{\operatorname{\sf Hdg}}}

\newcommand{\rk}{\operatorname{rk}}

\newcommand{\BK}{\operatorname{BK}}

\newcounter{Mycounter}[section]
\newcounter{lemma}[section]
\newcounter{claim}[section]
\newcounter{sublemma}[section]
\newcounter{corollary}[section]
\renewcommand{\thecorollary}{{Corollary \thesection.\arabic{corollary}}}
\newcommand{\corollary}{%
    \setcounter{corollary}{\value{Mycounter}}
    \refstepcounter{corollary}
    \stepcounter{Mycounter}
    {\noindent \bf \thecorollary:\ }}

\newcounter{theorem}[section]
\renewcommand{\thetheorem}{{Theorem \thesection.\arabic{theorem}}}
\newcommand{\theorem}{%
    \setcounter{theorem}{\value{Mycounter}}
    \refstepcounter{theorem}
    \stepcounter{Mycounter}
    {\noindent \bf \thetheorem:\ }}

\newcounter{conjecture}[section]
\newcounter{proposition}[section]
\renewcommand{\theproposition}
      {{Proposition \thesection.\arabic{proposition}}}
\newcommand{\proposition}{%
    \setcounter{proposition}{\value{Mycounter}}
    \refstepcounter{proposition}
    \stepcounter{Mycounter}
    {\noindent \bf \theproposition:\ }}

\newcounter{definition}[section]
\renewcommand{\thedefinition}
      {{Definition~\thesection.\arabic{definition}}}
\newcommand{\definition}{%
    \setcounter{definition}{\value{Mycounter}}
    \refstepcounter{definition}
    \stepcounter{Mycounter}
    {\noindent \bf \thedefinition:\ }}

\newcounter{example}[section]
\newcounter{remark}[section]
\renewcommand{\theremark}{{Remark \thesection.\arabic{remark}}}
\newcommand{\remark}{%
    \setcounter{remark}{\value{Mycounter}}
    \refstepcounter{remark}
    \stepcounter{Mycounter}
    {\noindent \bf \theremark:\ }}

\newcounter{problem}[section]
\newcounter{question}[section]
\makeatletter

\@addtoreset{equation}{section} \@addtoreset{footnote}{section}
\makeatother

\def\blacksquare{\hbox{\vrule width 5pt height 5pt depth 0pt}}
\def\endproof{\blacksquare}

\begin{document}

\begin{center}
{\LARGE\bf
Hyperbolic geometry of the ample cone\\[1mm] of a hyperk\"ahler manifold\\[1mm]
}

Ekaterina Amerik\footnote{Partially supported by 
RScF grant, project 14-21-00053, 11.08.14.}, Misha Verbitsky\footnote{Partially supported by 
RScF grant, project 14-21-00053, 11.08.14.


}

\end{center}

{\small \hspace{0.15\linewidth}
\begin{minipage}[t]{0.7\linewidth}
{\bf Abstract} \\
Let $M$ be a compact
hyperk\"ahler manifold with maximal holonomy
(IHS). The group $H^2(M, \R)$
is equipped with a quadratic form of signature $(3, b_2-3)$,
called Bogomolov-Beauville-Fujiki (BBF) form.
This form restricted to the rational Hodge lattice $H^{1,1}(M,\Q)$,
has signature $(1,k)$. This gives a hyperbolic Riemannian
metric on the projectivisation of the positive cone in $H^{1,1}(M,\Q)$,
denoted by $H$. Torelli theorem implies that the Hodge monodromy group $\Gamma$ acts on $H$ with finite covolume, giving
a hyperbolic orbifold $X=H/\Gamma$. We show that there are finitely many 
geodesic hypersurfaces which cut $X$ 
into finitely many polyhedral pieces in such a way that each of these pieces is 
isometric to a quotient $P(M')/\Aut(M')$, where $P(M')$
is the projectivization of the ample cone of a birational model $M'$ of $M$, 
and $\Aut(M')$ the group of its holomorphic automorphisms.
This is used to prove the existence of nef isotropic line bundles 
on a hyperk\"ahler birational model of a simple hyperk\"ahler manifold of 
Picard number at least 5, and also illustrates the fact that 
an IHS manifold 
has only finitely many birational models up to isomorphism (cf. \cite{MY}).
\end{minipage}
}

\tableofcontents


\section{Introduction}


Let $M$ be an irreducible holomorphically symplectic manifold, that is,
a simply-connected compact K\"ahler manifold with $H^{2,0}(M)=\C\Omega$
where $\Omega$ is nowhere degenerate. In dimension two, such manifolds are 
K3 surfaces; in higher dimension $2n,\ n>1,$ one knows, up to deformation, 
two infinite series of such manifolds, namely the punctual Hilbert schemes
of K3 surfaces and the generalized Kummer varieties, and two sporadic 
examples constructed by O'Grady. Though considerable effort has been made
to construct other examples, none is known at present, and the classification
problem for irreducible holomorphic symplectic manifolds (IHSM) looks equally 
out of reach.  

One of the main features of an IHSM $M$ is the existence of an integral 
quadratic form $q$ on the second cohomology $H^2(M,\Z)$, the 
{\bf Beauville-Bogomolov-Fujiki form} (BBF) form. It generalizes the intersection form
on a surface; in particular its signature is $(3, b_2-3)$, and the
signature of its restriction to $H^{1,1}_\R(M)$ is $(1, b_2-3)$. The cone
$\{x\in H^{1,1}_\R(M)|x^2>0\}$ thus has two connected components; we call
the {\it positive cone} $\Pos(X)$ the one which contains the K\"ahler classes.
The BBF form is in fact of topological origin: by a formula due to Fujiki,
$q(\alpha)^n$ is proportional to $\alpha^{2n}$ with a positive coefficient
 depending only on $M$. 
      
To understand better the geometry of an IHSM, it can be useful to fiber it, 
whenever possible,
over
a lower-dimensional variety. Note that by a result of Matsushita, the fibers
are always Lagrangian (in particular, $n$-dimensional, where $2n=\dim_\C M$),
and the general fiber is a torus. Such fibrations
are important for the classification-related problems, and
one can also hope to get some interesting geometry from
their degenerate fibers (for instance, use
them to construct rational curves on $M$).

Note that a fibration of $M$ is necessarily given by a linear system $|L|$ 
where $|L|$ is a nef line bundle with $q(L)=0$. Conjecturally, any such
bundle is semiample, that is, for large $m$ the linear system $|L^{\otimes m}|$ is base-point-free and thus gives a desired fibration.

It is therefore important to understand which irreducible holomorphic
symplectic varieties carry nef line bundles of square zero. By Meyer's theorem (see for example \cite{Se}),
$M$ has an integral $(1,1)$-class of square zero as soon as the Picard number
$\rho(M)$ is at least five. By definition, such a class is nef when it is 
in the closure of the K\"ahler cone $\Kah(M)\subset \Pos(M)$. The question is 
thus to understand whether one
can find an isotropic integral $(1,1)$-class in the closure of the K\"ahler 
cone. 

For projective K3 surfaces, this is easy and has been done in \cite{Psh-Sh}. 
Indeed $\Kah(M)\subset \Pos(M)$ is cut out by the orthogonal hyperplanes to 
$(-2)$-classes, since a positive $(1,1)$-class is K\"ahler if and only if
it restricts positively on all $(-2)$-curves, and $(-2)$-classes on a K3 surface are $\pm$-effective by Riemann-Roch. 
Let $x$ be an 
isotropic integral $(1,1)$-class and suppose that $x\not\in \overline{\Kah(M)}$,
that is, there is a $(-2)$-curve $p$ with $\langle x, p\rangle <0$.
Fix an ample integral $(1,1)$-class $h$. Then the image of $x$ under the 
reflection in $p^{\bot}$,  $x'=x+\langle x,p\rangle p$, satisfies 
$\langle x', h\rangle < \langle x, h\rangle$. Therefore the image of $x$ 
under successive reflexions in such $p$'s becomes nef after finitely many 
steps. The non-projective case is even easier, since an isotropic line bundle 
must then be in the kernel of the Neron-Severi lattice and so has zero 
intersection with every curve, in particular, it is nef.

Trying to apply the same argument to higher-dimensional IHSM we see that
the existence of an isotropic line bundle yields and isotropic element
in the closure of the {\bf birational K\"ahler cone}
$\BK(M)$. By definition, $\BK(M)$  is a
union of inverse images of the K\"ahler cone on all IHSM birational models
of $M$, and its closure is cut out in $\Pos(X)$ by the Beauville-Bogomolov 
orthogonals to the classes of the prime uniruled 
exceptional divisors (\cite{_Boucksom_}). One knows that 
the reflections
in those hyperplanes are integral (\cite{_Markman:reflections_}); in particular the divisors have bounded 
squares and the ``reflections argument'' above applies with obvious modifications.

A priori, the closure of $\BK(M)$ may strictly contain the union of the 
closures of the inverse images of the K\"ahler cones of all birational models,
so an additional argument is required to conclude that there is an isotropic
nef class on some birational model of $M$. One way to deal with this is
explained in the paper \cite{MZ}: the termination of flops on an IHSM implies
that any element of the closure of $\BK(M)$ does indeed become nef on some 
birational model. These observations, though, require the use of rather heavy 
machinery of the Minimal Model Program (MMP) which are in principle valid on all
varieties (though the termination of flops itself remains unproven in general).

The purpose of the present note is to give another proof of the existence
of nef isotropic classes, which does
not rely on the MMP. Instead it relies on the ``cone conjecture'' which was 
established
in \cite{_AV:Mor_Kaw_} using completely different methods, namely ergodic theory and
hyperbolic geometry. We find the hyperbolic geometry picture which appears in 
our proof particularly appealing, and believe that it might provide an 
alternative, perhaps sometimes more efficient,
approach to birational geometry in the particular 
case of the irreducible holomorphic 
symplectic manifolds.

The main advantage of the present construction is
its geometric interpretation. 
The BBF quadratic form, restricted to the rational Hodge lattice $H^{1,1}(M,\Q)$,
has signature $(1,k)$ (unless $M$ is non-algebraic,
in which case our results are tautologies). 
This gives a hyperbolic Riemannian
metric on the projectivisation of the positive cone in $H^{1,1}(M,\Q)$,
denoted by $H$. Torelli theorem implies that the group $\Gamma^\Hdg$
of Hodge monodromy acts on $H$ with finite covolume, giving
a hyperbolic orbifold $X=H/\Gamma^\Hdg$. Using Selberg lemma, one easily reduces to the case when $X$ is a manifold. We prove that $X$ is cut
into finitely many polyhedral pieces by finitely many
geodesic hypersurfaces in such a way that each of these pieces is 
isometric to a quotient $\Amp(M')/\Aut(M')$, where $\Amp(M')$
is the projectivization of the ample cone of a birational model of $M$, 
and $\Aut(M')$ the group of holomorphic automorphisms.

In this interpretation, equivalence classes of birational
models are in bijective correspondence with these 
polyhedral pieces $H_i$, and the isotropic nef line bundles
correspond to the cusp points of these $H_i$.
Existence of cusp points is implied by Meyer's theorem,
and finiteness of $H_i$ 
by our results on the cone conjecture from \cite{_AV:Mor_Kaw_} (Section \ref{_Hyper_geo_Section_}). 
Finally, the geometric finiteness results from hyperbolic geometry imply the finiteness of the isotropic nef line
bundles up to automorphisms.

\section{Hyperk\"ahler manifolds: basic results}

In this section, we recall the definitions and basic
properties of hyperk\"ahler manifolds and MBM classes.

\subsection{Hyperk\"ahler manifolds}\label{BBF}

\definition
A {\bf hyperk\"ahler manifold} $M$, that is, a compact K\"ahler holomorphically symplectic manifold,
is called {\bf simple} (alternatively, {\bf irreducible holomorphically symplectic (IHSM)}), if 
$\pi_1(M)=0$ and $H^{2,0}(M)=\C$.

\hfill



This definition is motivated by the following theorem
of Bogomolov. 

\hfill

\theorem \label{_Bogo_deco_Theorem_}
(\cite{_Bogomolov:decompo_})
Any hyperk\"ahler manifold admits a finite covering
which is a product of a torus and several 
simple hyperk\"ahler manifolds.
\endproof

\hfill

The second cohomology $H^2(M,\Z)$ of a simple 
hyperk\"ahler manifold $M$ carries a
primitive integral quadratic form $q$, 
called {\bf the Bogomolov-Beauville-Fujiki form}. It generalizes the intersection product on a K3 surface:
its signature is $(3,b_2-3)$ on $H^2(M,\R)$ and $(1,b_2-3)$ on $H^{1,1}_{\R}(M)$.  It was first
defined in \cite{_Bogomolov:defo_} and 
\cite{_Beauville_},
but it is easiest to describe it using the
Fujiki theorem, proved in \cite{_Fujiki:HK_}.

\hfill

\theorem\label{_Fujiki_Theorem_}
(Fujiki)
Let $M$ be a simple hyperk\"ahler manifold,
$\eta\in H^2(M)$, and $n=\frac 1 2 \dim M$. 
Then $\int_M \eta^{2n}=c q(\eta,\eta)^n$,
where $q$ is a primitive integer quadratic form on $H^2(M,\Z)$, and 
$c>0$ is a rational number. \endproof

\hfill

%

\definition
Let $M$ be a hyperk\"ahler manifold. The 
{\bf monodromy group} of $M$ is a subgroup of $GL(H^2(M,\Z))$
generated by the monodromy transforms for all Gauss-Manin local systems.

\hfill

It is often enlightening to consider this group in terms of the
mapping class group action. We briefly recall this description. 


\hfill

The {\bf Teichm\"uller space} $\Teich$ is the quotient $\Comp(M)/\Diff_0(M)$, where $\Comp(M)$
denotes the space of all complex structures of K\"ahler type on $M$ and $\Diff_0(M)$ is the group of isotopies.
It follows from a result of Huybrechts (see \cite{_Huybrechts:finiteness_}) that for an IHSM $M$, $\Teich$ has
only finitely many connected components. Let $\Teich_M$ denote the one containing our given complex structure.
Consider the subgroup of the mapping class group $\Diff(M)/\Diff_0(M)$ fixing $\Teich_M$.

\hfill

\definition\label{monodr} The {\bf monodromy group}
$\Gamma$ is the image of this subgroup in $O(H^2(M,\Z), q)$. The {\bf Hodge monodromy group} $\Gamma^{\Hdg}$
is the subgroup of $\Gamma$ preserving the Hodge decomposition.







%


\hfill

\theorem\label{arithmetic} (\cite{_V:Torelli_}, Theorem 3.5) 
The monodromy group is a finite index subgroup in $O(H^2(M, \Z), q)$
(and the Hodge monodromy is therefore an arithmetic
subgroup of the orthogonal group of the 
Picard lattice).

\subsection{MBM classes}



\definition
A cohomology class $\eta\in H^2(M, \R)$ is called 
{\bf positive} if $q(\eta,\eta)>0$, and 
{\bf negative} if $q(\eta,\eta)<0$. The {\bf positive cone} $\Pos(M)\in H^{1,1}_{\R}(M)$
is that one of the two connected components of the set 
of positive classes on $M$ which
contains the K\"ahler classes.

\hfill

Recall e.g. from \cite{_Markman:survey_} that the positive cone is decomposed into the union of 
{\bf birational K\"ahler chambers}, which are monodromy transforms of the 
{\bf birational K\"ahler cone} $\BK(M)$. The birational K\"ahler cone is,
by definition, the union of pullbacks of the K\"ahler cones $\Kah(M')$ where $M'$
denote a hyperk\"ahler birational model of $M$ (the ``K\"ahler chambers'').  
The {\bf faces}\footnote{A face of a convex cone in a vector space $V$
is the intersection of its boundary and a hyperplane which 
has non-empty interior in the hyperplane.}  
of these chambers are supported on the hyperplanes
orthogonal to the classes of prime
uniruled divisors of negative square on $M$.


\hfill

The {\bf MBM classes} are defined as those classes whose orthogonal hyperplanes support
faces of the K\"ahler chambers.

\hfill

\definition\label{mbm}
A negative integral cohomology class $z$ of type $(1,1)$
is called {\bf monodromy birationally minimal} (MBM)
if for some isometry $\gamma\in O(H^2(M,\Z))$ 
belonging to the monodromy group,
 $\gamma(z)^{\bot}\subset H^{1,1}_{\R}(M)$ contains a face 
of the K\"ahler cone of one of birational
models $M'$ of $M$.

\hfill

Geometrically, the MBM classes are characterized among negative integral 
$(1,1)$-classes, as those which are, up to a scalar multiple, represented 
by minimal rational curves on deformations of $M$ under the identification of $H_2(M,\Q)$ with $H^2(M,\Q)$ given by the BBF form (\cite{_AV:MBM_}, \cite{AV3}, \cite{KLM}). 

\hfill

The following theorems summarize the main results about MBM classes from \cite{_AV:MBM_}. 

\hfill

\theorem\label{defo-inv}(\cite{_AV:MBM_}, Corollary 5.13)
An MBM class $z\in H^{1,1}(M)$ is also MBM on any deformation $M'$ of $M$ where $z$ remains of
type $(1,1)$.

\hfill

\theorem\label{kahler-cone}(\cite{_AV:MBM_}, Theorem 6.2)
The K\"ahler cone of $M$ is 
a connected component of $\Pos(M)\backslash \cup_{z\in S} z^\bot$,
where $S$ is the set of MBM classes on $M$.

\hfill

In what follows, we shall also consider the positive cone in the algebraic part 
$NS(M)\otimes \R$ of $H^{1,1}_{\R}(M)$, denoted by $\Pos_{\Q}(M)$. Here and further on, $NS(M)$ stands
for N\'eron-Severi group of $M$.

\hfill

\definition The {\bf ample chambers} are the connected components of 
$\Pos_{\Q}(M)\backslash \cup_{z\in S} z^\bot$ where $S$ is the set of MBM classes on $M$. 

\hfill

One of the ample chambers is, obviously, the ample cone of $M$, hence the name.

\hfill

In the same way, one defines {\bf birationally ample} or {\bf movable} chambers as the connected components
of the complement to the union of orthogonals to the classes of uniruled divisors and their monodromy transforms,
cf. \cite{_Markman:survey_}, section 6. These are also described as intersections of the biratonal K\"ahler
chambers with $NS(M)\otimes \R$.

\hfill

\remark Because of the deformation-invariance property of
 MBM classes, it is natural to introduce this notion on $H^2(M,\Z)$ rather than on $(1,1)$-classes:
we call $z\in H^2(M,\Z)$ an MBM class
as soon as it is MBM in those complex structures where it is of 
type $(1,1)$.





\subsection{Morrison-Kawamata cone conjecture}

The following theorem has been proved in \cite{_AV:Mor_Kaw_}.

\hfill

\theorem (\cite{_AV:Mor_Kaw_}) Suppose that the Picard number $\rho(M)>3$. Then
the Hodge monodromy group has only finitely many orbits on the set of MBM classes
of type $(1,1)$ on $M$.
\endproof

\hfill

Since the Hodge monodromy group acts by isometries, it follows that the primitive
MBM classes have bounded square (using the deformation argument, one easily extends
this last statement from the case of $\rho(M)>3$ to that of $b_2(M)\neq 5$, but we
shall not need this here). In \cite{_AV:MBM_} we have seen that this implies
some apriori stronger statements on the Hodge monodromy action.

\hfill

\corollary \label{finiteness-faces} The Hodge monodromy group has only finitely
many orbits on the set of faces of the K\"ahler chambers, as well as on the set of the K\"ahler
chambers themselves.

\hfill

For reader's convenience, let us briefly sketch the proof (for details, see sections 3 and 6 of \cite{_AV:MBM_}). It consists in remarking 
that a face of 
a chamber is given by a flag
$P_s\supset P_{s-1}\supset \dots \supset P_1$ where $P_s$ is the supporting hyperplane
(of dimension $s=h^{1,1}-1$), $P_{s-1}$ supports a face of our face, etc., and for each
$P_i$ an {\bf orientation} (``pointing inwards the chamber'') is fixed. One deduces
from the boundedness of the square of primitive MBM classes that possible $P_{s-1}$
are as well given inside $P_s$ by orthogonals to integral vectors of bounded square,
and it follows that the stabilizer of $P_s$ in $\Gamma^{\Hdg}$ acts with finitely many
orbits on those vectors; continuing in this way one eventually gets the statement.

\hfill

By Markman's version of the Torelli theorem \cite{_Markman:survey_}, an element of $\Gamma^{\Hdg}$ preserving
the K\"ahler cone actually comes from an automorphism of $M$. Thus an immediate consequence
is the following K\"ahler version of the Morrison-Kawamata cone conjecture.

\hfill

\corollary (\cite{_AV:Mor_Kaw_}) $\Aut(M)$ has only finitely many orbits on the set of
faces of the K\"ahler cone.

\hfill

\remark\label{faces-ample} As the faces of the ample cone are likewise given by the orthogonals to MBM classes, 
but in $\Pos_{\Q}(M)$ rather than in $\Pos(M)$, one concludes that the same must be true for the
ample cone.


\section{Hyperbolic geometry and the K\"ahler cone}
\label{_Hyper_geo_Section_}


\subsection{Kleinian groups and hyperbolic manifolds}

\definition
A {\bf Kleinian group} is a discrete subgroup of isometries of the hyperbolic
space $\H^n$.

\hfill

One way to view $\H^n$ is as a projectivization of the positive cone $\P V^+$
of a quadratic form $q$ of signature $(1,n)$ on a real vector space $V$. The Kleinian
groups are thus discrete subgroups of $SO(1,n)$. One calls such a subgroup a {\bf lattice}
if its covolume is finite.

\hfill

\definition
An {\bf arithmetic subgroup} of an algebraic group
$G$ defined over the integers is a subgroup commensurable with $G_\Z$.

\hfill

\remark\label{BHC} From Borel and Harish-Chandra theorem (see
\cite{bo-hch}) it follows that when $q$ is integral,
any arithmetic subgroup of $SO(1,n)$ is a lattice for $n\geq 2$.

\hfill

\definition A {\bf complete hyperbolic orbifold} is a quotient of
the hyperbolic space 
by a Kleinian group. A {\bf complete hyperbolic manifold} is a quotient
of the hyperbolic space by a Kleinian group acting freely.

\hfill

\remark One defines a hyperbolic manifold as 
a manifold of constant negative bisectional curvature.
When complete, such a manifold is uniformized 
by the hyperbolic space (\cite{_Thurston:thick-thin_}).

\hfill

The following proposition is well-known.

\hfill

\proposition\label{mfd} Any complete hyperbolic 
orbifold has a finite covering which is a
complete hyperbolic manifold (in other words, 
any Kleinian group has a finite index
subgroup acting freely).

\hfill

{\bf Proof:} Let $\Gamma$ be a Kleinian group. Notice first that all stabilizers
for the action of $\Gamma$ on $\P V^+$ are finite, since these are identified to
discrete subgroups of a compact group $SO(n)$. Now by Selberg lemma  $\Gamma$ has
a finite index subgroup without torsion which must therefore act freely.

\hfill

\remark If $M$ is an IHSM, the group of Hodge monodromy $\Gamma^{\Hdg}$
is an arithmetic lattice in $SO(H^{1,1}(M,\Q))$ when $\rk
H^{1,1}(M,\Q)\geq 3$. The hyperbolic manifold 
${\Bbb P}(H^{1,1}(M,\Q)\otimes_\Q\R)^+/\Gamma^\Hdg$ 
has finite volume by Borel and Harish-Chandra theorem.

\subsection{The cone conjecture and hyperbolic geometry}

Recall that {\bf the rational positive cone}
$\Pos_\Q(M)$ of a projective hyperk\"ahler manifold $M$
is one of two connected components of the set of positive
vectors in $NS(M)\otimes \R$.

Replacing $\Gamma^{\Hdg}$ by a finite index subgroup if necessary, we may
assume that the quotient $\P \Pos_\Q(M)/ \Gamma^{\Hdg}$ is a complete hyperbolic 
manifold which we shall denote by $H$.

By Borel and Harish-Chandra theorem (see \ref{BHC}), $H$ is of finite volume as soon as
the Picard number of $M$ is at least three.

Let  $S=\{s_i\}$ be the set of MBM classes of type $(1,1)$ on $M$.
The following is a translation of the Morrison-Kawamata cone conjecture into the setting
of hyperbolic geometry.

\hfill

\theorem The images of the hyperplanes $s_i^{\bot}$,
$s_i\in S$, cut $H= {\Bbb P}\Pos_\Q(M)/\Gamma^{Hdg}$ into finitely many
pieces. One of those pieces is the image of the ample cone (up to a finite covering, this is the 
quotient of the ample cone by $\Aut(M)$) and the others are the images of 
ample cones of birational models of $M$. The closure of each one is a hyperbolic manifold with
boundary consisting of finitely many geodesic pieces.

\hfill

{\bf Proof:} According to \ref{finiteness-faces}, up to the action of $\Gamma^{\Hdg}$  
there are finitely many faces of ample chambers. Each face is a connected component 
of the complement to $\cup_{j\neq i} s_j^{\bot}$ in $s_i^{\bot}$ for some $i$. It is clear that the images 
of the faces do not intersect hence, being finitely many, cut $H$ into finitely many pieces which are images of the ample 
chambers. We have already mentioned that an element of $\Gamma^{\Hdg}$ preserving the K\"ahler cone
is induced by an automorphism. Finally, the whole $H$ is covered by the birational ample cone (since the other 
birational ample chambers are its monodromy transforms) and thus each part of $H$ obtained in this
way comes from an ample chamber. \endproof

\hfill






%








Let us also mention that the same arguments also prove the following result (cf. \cite{MY}).

\hfill

\corollary
There are only finitely many non-isomorphic birational
models of $M$.

\hfill

{\bf Proof:}
Indeed, the K\"ahler (or ample) chambers in the same  $\Gamma^{\Hdg}$-orbit correspond to isomorphic 
birational models, since one can view the action of  $\Gamma^{\Hdg}$ as the change of the marking
(recall that a {\bf marking} is a choice of an isometry of $H^2(M,\Z)$ with a fixed lattice $\Lambda$
and that there exists a coarse moduli space of marked IHSM which in many works (e.g. \cite{_Huybrechts:basic_}) plays
the same role as the Teichm\"uller space in others).
\endproof



\section{Cusps and nef parabolic classes}

\definition
{\bf A horosphere} on a hyperbolic space is 
a sphere which is everywhere orthogonal to a pencil
of geodesics passing through one point at infinity, and 
{\bf a horoball} is a ball bounded by a horosphere.
{\bf A cusp point} for an $n$-dimensional hyperbolic manifold ${\Bbb
  H}/\Gamma$ is a point on the boundary $\6{\Bbb H}$ 
such that its stabilizer in $\Gamma$ 
contains a free abelian group of rank $n-1$.
Such subgroups are called {\bf maximal parabolic}.
For any point $p\in \6{\Bbb H}$ 
stabilized by $\Gamma_0\subset \Gamma$,
and any horosphere $S$ tangent to the boundary in $p$,
$\Gamma_0$ acts on $S$ by isometries. 
In such a situation, $p$ is a cusp point 
if and only if $(S\backslash p)/\Gamma_0$
is compact.

\hfill

A cusp point $p$ yields a {\bf cusp} 
in the quotient ${\Bbb H}/\Gamma$, that is,
a geometric end of ${\Bbb H}/\Gamma$ 
of the form $B/\Z^{n-1}$, where 
$B\subset \H$ is a horoball tangent
to the boundary at $p$.

The following theorem describes the geometry of finite
volume complete hyperbolic manifolds more precisely.

\hfill

\theorem (Thick-thin decomposition)
 \\
Any $n$-dimensional complete hyperbolic 
manifold of finite volume can be represented as
a union of a ``thick part'', which is a compact manifold
with a boundary, and a ``thin part'', which is
a finite union of quotients of form $B/\Z^{n-1}$, where $B$ is a horoball
tangent to the boundary at a cusp point, and $\Z^{n-1} =\St_\Gamma(B)$.

{\bf Proof:} See \cite[Section 5.10]{_Thurston:thick-thin_}
or \cite[page 491]{_Kapovich:Kleinian_}). \endproof 

\hfill

\theorem\label{_cusps_rat_points_Theorem_}
Let ${\Bbb H}/\Gamma$ be a hyperbolic manifold,
where $\Gamma$ is an arithmetic subgroup of $SO(1,n)$.
Then the cusps of ${\Bbb H}/\Gamma$ are in (1,1)-correspondence with
$Z/\Gamma$, where $Z$ is the set of rational lines
$l$ such that $l^2=0$.

\hfill

{\bf Proof:}
By definition of cusp points, the
cusps of ${\Bbb H}/\Gamma$ are in 1 to 1 
correspondence with $\Gamma$-conjugacy classes of 
maximal parabolic subgroups of $\Gamma$
(see \cite{_Kapovich:Kleinian_}). Each such
subgroup is uniquely determined by the
unique point it fixes on the boundary of ${\Bbb H}$.
\endproof

\hfill

The main result of this paper is the following theorem.

\hfill

\theorem
Let $M$ be a hyperk\"ahler manifold with Picard number at least 5. Then $M$ has a birational model
admitting an integral nef $(1,1)$-class $\eta$ with $\eta^2 =0$. Moreover each birational model
contains only finitely many such classes up to automorphism.

\hfill

{\bf Proof:} By Meyer's theorem (see for example \cite{Se}), there exists
$\eta\in NS(M)$ with $\eta^2=0$. 
By \ref{_cusps_rat_points_Theorem_}, the hyperbolic manifold
$H:={\Bbb P}\Pos_\Q(M)/\Gamma^{\Hdg}$ then has cusps, and, being of finite volume, only finitely many of them. 
Recall that $H$ is decomposed into finitely many pieces, and each of those pieces is the image of 
the ample cone of a birational model of $M$ in $\Pos_\Q(M)$. Therefore a lifting of each cusp to the boundary of 
${\Bbb P}\Pos_\Q(M)$ gives a BBF-isotropic nef line bundle on a birational model $M'$ of $M$ (or, more precisely, the
whole line such a bundle generates in $NS(M)\otimes \R$). Finally, the number of $\Aut(M')$-orbits of such 
classes is finite, being exactly the number of cusps  in the piece of $H$ corresponding to $M'$: indeed this piece is 
just the quotient of the ample cone of $M'$ by its
stabilizer which is identified with $\Aut(M')$.
\endproof

\hfill

{\bf Acknowledgements:}
We are grateful to S. Cantat, M. Kapovich and V. Gritsenko
 for interesting discussions and advice.

{

\noindent {\sc Ekaterina Amerik\\
{\sc Laboratory of Algebraic Geometry,\\
National Research University HSE,\\
Department of Mathematics, 7 Vavilova Str. Moscow, Russia,}\\
\tt  Ekaterina.Amerik@gmail.com}, also: \\
{\sc Universit\'e Paris-11,\\
Laboratoire de Math\'ematiques,\\
Campus d'Orsay, B\^atiment 425, 91405 Orsay, France}

\hfill

\noindent {\sc Misha Verbitsky\\
{\sc Laboratory of Algebraic Geometry,\\
National Research University HSE,\\
Department of Mathematics, 7 Vavilova Str. Moscow, Russia,}\\
\tt  verbit@mccme.ru}, also: \\
{\sc Universit\'e Libre de Bruxelles, CP 218,\\
Bd du Triomphe, 1050 Brussels, Belgium}
 }

\end{document}